\documentclass[12pt, A4 paper]{article}
\usepackage{amssymb}
\usepackage{amsmath}
\usepackage{graphicx}
\usepackage[authoryear]{natbib}
\usepackage{times,epsfig,makeidx,amsfonts,amsmath,dcolumn,enumerate,multirow,graphicx,amssymb,colortbl}
\usepackage[margin=1in]{geometry}

\newtheorem{theorem}{Theorem}
\newtheorem{corollary}{Corollary}

\newenvironment{proof}[1][Proof]{\noindent\textbf{#1.} }{\ \rule{0.5em}{0.5em}}

\newcommand{\qed}{\hspace{\stretch{3}}$\square$\\[1.8ex]}

\begin{document}
\title{A note on the infinite divisibility of a class of transformations of normal variables}
\author{{\large A. Murillo-Salas}\footnote{{\sc Departamento de Matem\'{a}ticas, Universidad de Guanajuato, Jalisco S/N, Mineral de Valenciana, Guanajuato, Gto.\ C.P.\ 36240,
M\'exico.} E-mail: amurillos@ugto.mx} \,{\large and}  {\large F. J. Rubio}\footnote{{\sc University of Warwick, Department of Statistics, Coventry, CV4 7AL.} E-mail:  F.J.Rubio@warwick.ac.uk}}
\maketitle

\begin{abstract}
This note examines the infinite divisibility of density-based transformations of normal random variables. We characterize a class of density-based transformations of normal variables which produces non-infinitely divisible distributions. We relate our result with some known skewing mechanisms.
\end{abstract}



%------------------------------------------------------------------------------------------------------------------------------------------------------------------------------------------------------

\section{Introduction}

In recent years, the use of skewed distributions has attracted the attention of statisticians both from the applied and theoretical points of view. Several skewing mechanisms have been proposed to obtain skewed distributions by transforming symmetric ones \citep{AValle05,Jones04,MO97,Wang04}.
Thus, it is of interest {\it per se} to analyze theoretical properties of the transformed distributions. Here, we present a
characterization (see Theorem \ref{MainTh} below) for a class of density-based transformations that produces non-infinitely divisible
 distributions when applied to the normal distribution. This characterization can be easily related with some skewing mechanisms. It is worth mentioning that our result has implications in statistical modeling because it rules out the use of several skew-normal distributions in
 models defined in terms of infinitely divisible distributions (cf. \cite{Steu79}).

In Section $\ref{S-M}$, a general representation of density-based transformations proposed in \cite{FS06} and its relationship with four skewing
mechanisms is presented. Using this representation, in Section $\ref{Div}$ we offer a sufficient condition on such density-based transformations
which destroys the property of infinite divisibility of normal distributions. Our results partially extend the work of \cite{Dom07} and \cite{Koz08}. We specialize our characterization to some skewing mechanisms  which are of interest in the statistical literature.

%------------------------------------------------------------------------------------------------------------------------------------------------------------------------------------------------------

\section{Skewing mechanisms}\label{S-M}

Let $S$ and $F$ be two distribution functions on the real line and let $P$ be a distribution function on $(0,1)$,
densities (which are assumed to exist) will be denoted by the corresponding lowercase letters. \cite{FS06} show that for any pair of absolutely continuous distributions $S$ and $F$ with support on ${\mathbb R}$, there exists a distribution $P$ such that $S=P\circ F$. This implies that the transformation from a random variable $X$ with distribution $F$ to a random variable $Y$ with distribution $S$ can be represented as a density-based transformation as follows

\begin{eqnarray}\label{representation}
s(y\vert F,P) = f(y)p[F(y)],
\end{eqnarray}

\noindent $y \in {\mathbb R}$. If the distribution $F$ is symmetric and $S$ is asymmetric, then $S$ is said to be an asymmetric version of $F$ generated
by the skewing mechanism $P$ \citep{FS06}. Several skewing mechanisms have been proposed in the statistical literature. Next, we provide the relationship
of four known skewing mechanisms with representation $(\ref{representation})$:

\begin{enumerate}[(i)]
\item Skew-symmetric construction \citep{Wang04}. Such transformation is defined as follows
\begin{eqnarray*}
s(y \vert F,P) &=& 2f(y)\pi(y),\\
p(x) &=& 2\pi(F^{-1}(x)),
\end{eqnarray*}
where $\pi$ is a function that satisfies $0\leq \pi(y) \leq 1$, $\pi(-y)=1-\pi(y)$ and $x \in [0,1]$.

Particular cases of this transformation are the Hidden Truncation skewing mechanism \citep{Arnold00} and Azzalini's skew-normal \citep{Azzalini85}.

\item Order Statistics \citep{Jones04}. This transformation introduces two new parameters $\psi_1>0$ and $\psi_2>0$ as follows

\begin{eqnarray*}
s(y \vert F,P) &=&  [\beta( \psi_1,\psi_2)]^{-1}f(y)F(y)^{\psi_1-1}(1-F(y))^{\psi_2-1},\\
p(x \vert \psi_1,\psi_2) &=& [\beta( \psi_1,\psi_2)]^{-1}x^{\psi_1-1}(1-x)^{\psi_2-1},
\end{eqnarray*}

\noindent where $\beta$ denotes the beta function and $x \in [0,1]$.

\item  Marshall-Olkin transformation \citep{MO97}. Such transformation is given through a positive parameter $\gamma$ as follows

\begin{eqnarray*}
s(y\vert F,P) &=& \frac{\gamma f(y)}{[F(y)+\gamma(1-F(y))]^2},\\
p(x \vert \gamma) &=& \frac{\gamma}{[x + \gamma(1-x)]^2},
\end{eqnarray*}

\noindent where $x \in [0,1]$.

\item Two-piece distributions \citep{AValle05}. This transformation consists of scaling by different factors, $a(\gamma)$ and $b(\gamma)$,
 either side of the symmetry point of a unimodal density $f$. If $f$ is symmetric about $0$, this transformation is given by

\begin{eqnarray*}
s(y\vert F, P) &=& \frac{2}{a(\gamma)+b(\gamma)}\left[ f\left(\frac{y}{b(\gamma)}\right)I(y<0)+f\left(\frac{y}{a(\gamma)}\right)I(y\geq 0)\right],\\
p(x \vert \gamma) &=& \frac{2}{a(\gamma)+b(\gamma)}\frac{f\left(\frac{F^{-1}(x)}{b(\gamma)}\right)I(x<0.5)+f\left(\frac{F^{-1}(x)}{a(\gamma)}\right)I(x\geq 0.5)}{f\left(F^{-1}(x)\right)},
\end{eqnarray*}
with $x \in [0,1]$, $a(\gamma)$ and $b(\gamma)$ are positive functions of the parameter $\gamma \in \Gamma$; where $\Gamma$ depends on the choice of the functions $\{a(\gamma),b(\gamma)\}$.

This class of transformations includes the Inverse Scale Factors presented in \cite{FS98} and the $\epsilon-$skew normal given in \cite{Mud00}.
\end{enumerate}

%------------------------------------------------------------------------------------------------------------------------------------------------------------------------------------------------------
\section{On infinite divisibility}\label{Div}

The goal of this section is twofold. Firstly, we consider $f$ to be normal in $(\ref{representation})$, then we derive a sufficient
condition on the density $p$ such that the skewed distribution $S$ is not infinitely divisible. Secondly, we relate this
result with the skewing mechanisms described in Section \ref{S-M}.

%----------------------------------------------------------------------------------------------------------------------------------------------------------------------------------------------------
%----------------------------------------------------------------------------------------------------------------------------------------------------------------------------------------------------

\begin{theorem}\label{MainTh}

Let $\Phi$ and $\phi$ be the distribution and density functions of a standard normal variable, respectively. Consider
 the transformation given in $(\ref{representation})$, with $p$ bounded {\sl a.e.}. Then, $S$ is
not infinitely divisible unless it is normal.
\end{theorem}
\proof Suppose that $p\leq M$, where $M>0$. Note that for $y>0$

\begin{eqnarray*}
S(-y \vert \Phi, P) &=& \int_{-\infty}^{-y}\phi(t)p[\Phi(t)]\,{d t} \leq M\int_{-\infty}^{-y}\phi(t)\,d t = M \Phi(-y),\\
1 - S(y \vert \Phi, P) &=&  \int_{y}^{\infty}\phi(t)p[\Phi(t)]\,d t \leq M\int_{y}^{\infty}\phi(t)\,d t = M[1-\Phi(y)].
\end{eqnarray*}
Then
\begin{eqnarray*}
S(-y \vert \Phi, P) + 1 - S(y \vert \Phi, P) \leq M[\Phi(-y) + 1 - \Phi(y)].
\end{eqnarray*}
Therefore for $y>1$
\begin{eqnarray*}
-\frac{\log[S(-y \vert \Phi, P) + 1 - S(y \vert \Phi, P)]}{y\log(y)} \geq -\frac{\log[M[\Phi(-y) + 1 - \Phi(y)]]}{y\log(y)}.
\end{eqnarray*}
Hence, taking limits on both sides we get
\begin{eqnarray*}
\limsup_{y\rightarrow\infty}-\frac{\log[S(-y \vert \Phi, P) + 1 - S(y \vert \Phi, P)]}{y\log(y)} \geq \limsup_{y\rightarrow\infty}-\frac{\log[M[\Phi(-y) + 1 - \Phi(y)]]}{y\log(y)},
\end{eqnarray*}
which together with the characterization of the normal distribution given in \cite{Steu04} Corollary $9.9$ implies the
result.\qed

%----------------------------------------------------------------------------------------------------------------------------------------------------------------------------------------------------
Theorem \ref{MainTh} has the following immediate

\begin{corollary} \label{IDSM}
Skew normals obtained by the skew-symmetric construction, the Marshall-Olkin transformation and the Order Statistics transformation for $\psi_1>1$ and $\psi_2>1$, are non-infinitely divisible.
\end{corollary}
\proof It is enough to note that, in each case, the corresponding density $p$ is bounded.\qed

\cite{Dom07} and \cite{Koz08} prove that, in particular, Azzalini's skew normal is not infinitely divisible. Theorem $\ref{MainTh}$ together with
Corollary $\ref{IDSM}$ extend this result to the family of skew-normal distributions obtained by the skew-symmetric construction \citep{Wang04}, from
which Azzalini's skew normal is a particular case. Moreover, the skewed normal distribution obtained by any skewing mechanism which satisfies the condition given in Theorem
 $\ref{MainTh}$ will lose the  infinite divisibility property.

%---------------------------------------------------------------------------------------------------------------------------------------------------------------------------------------------------------

Note that, for the skewing mechanism that produces two-piece distributions, the corresponding $p$ is not necessarily bounded.
Thus, an {\sl ad hoc} proof of the non-infinite divisibility of the skew-normals obtained with this sort of transformation is presented in the
following
\begin{theorem}
The two-piece skew normal is non-infinitely divisible unless $a(\gamma)=b(\gamma)$.
\end{theorem}
\proof \cite{Jones06} proves that the elements of the class of two-piece distributions are reparameterizations of each other. Therefore it is enough to
prove the result for the particular choice $\{a(\gamma),b(\gamma)\} = \{1-\gamma,1+\gamma\}$, $\gamma \in (-1,1)$, analyzed in \cite{AValle05}. Note
that $a(\gamma)=b(\gamma)$ if and only if $\gamma=0$, which corresponds to the symmetric normal which is infinitely divisible.  In addition, the density $s$ obtained with a particular $\gamma$, corresponds to reflecting the density $s$ with parameter $-\gamma$ around $0$. Hence, it is enough to
prove the result for $-1<\gamma<0$.

Note that
\begin{eqnarray*}
S(y\vert \Phi, P)=(1+\gamma)\Phi\left(\frac{y}{1+\gamma}\right)I(y<0) + \left[ -\gamma + (1-\gamma)\Phi\left(\frac{y}{1-\gamma}\right)\right]I(y\geq 0).
\end{eqnarray*}
Then, given that $-1<\gamma<0$ and for any $y>0$ we have that
\begin{eqnarray*}
S(-y\vert \Phi, P)+ 1 - S(y\vert \Phi, P) &=& (1+\gamma)\Phi\left(-\frac{y}{1+\gamma}\right) + 1 + \gamma - (1-\gamma)\Phi\left(\frac{y}{1-\gamma}\right)\\
&<& 2(1-\gamma)\Phi\left(-\frac{y}{1-\gamma}\right) < 4\Phi\left(-\frac{y}{2}\right).
\end{eqnarray*}
Then, for $y>1$ we have
\begin{eqnarray*}
-\frac{\log\left[S(-y\vert \Phi, P)+ 1 - S(y\vert \Phi, P)\right]}{y\log(y)}&>& -\frac{ \log\left[4\Phi\left(-\frac{y}{2}\right)\right]}{y\log(y)}\\
&=& -\frac{ \log\left[ 2 \left(\Phi\left(-\frac{y}{2}\right) + 1 - \Phi\left(\frac{y}{2}\right) \right) \right]}{y\log(y)}.
\end{eqnarray*}
The result follows by taking limits in both sides of this expression and using the characterization of the normal distribution
given by \cite{Steu04}, Corollary $9.9$.\qed

\section*{Acknowledgements} The authors are indebted to V\'ictor P\'erez-Abreu for several useful comments which significantly improved
the presentation of the work. F. J. Rubio acknowledges support from CONACYT-M{\'e}xico.

%---------------------------------------------------------------------------------------------------------------------------------------------------------------------------------------------------------------


\begin{thebibliography}{9}

\bibitem[Arellano-Valle et al.(2005)]{AValle05} Arellano-Valle, R. B., G{\'o}mez, H. W. and Quintana, F. A. (2005). Statistical inference for a general class of asymmetric distributions. {\sl Journal of Statistical Planning and Inference} 128: 427--443.

\bibitem[Arnold and Beaver(2000)]{Arnold00} Arnold, B. C. and Beaver, R. J. (2000). Hidden truncation models. {\sl Sankhya A} 62: 23--35.

\bibitem[Azzalini(1985)]{Azzalini85} Azzalini, A. (1985). A class of distributions which includes the normal ones. {\sl Scandinavian Journal of Statistics} 12: 171-178.

\bibitem[Dom\'{i}nguez-Molina and Rocha-Arteaga (2007)]{Dom07} Dom\'{i}nguez-Molina, J. A. and Rocha-Ortega, A. (2007). On the infinite divisibility of some skewed symmetric distributions. {\sl Statistics $\&$ Probability Letters} 77: 644--648.

\bibitem[Fern{\'a}ndez and Steel(1998)]{FS98} Fern{\'a}ndez, C. and Steel, M. F. J. (1998). On Bayesian modeling of fat tails and skewness. {\sl Journal of the American Statistical Association} 93: 359--371.

\bibitem[Ferreira and Steel(2006)]{FS06}  Ferreira, J. T. A. S. and Steel, M. F. J. (2006). A constructive representation of univariate skewed distributions. {\sl Journal of the American Statistical Association}, Theory and Methods 101: 823--829.

\bibitem[Jones(2004)]{Jones04} Jones, M. C. (2004). Families of distributions arising from distributions of order statistics (with discussion). {\sl Test} 13: 1-43.

\bibitem[Jones(2006)]{Jones06} Jones, M. C. (2006). A note on rescalings, reparametrizations and classes of distributions. {\sl Journal of Statistical Planning and Inference} 136: 3730-3733.

\bibitem[Kozubowski and Nolan(2008)]{Koz08} Kozubowski, T. J. and Nolan J. P. (2008). Infinite divisibility of skew Gaussian and Laplace laws.{\sl Statistics $\&$ Probability Letters} 78: 654--660.

\bibitem[Marshall and Olkin(1997)]{MO97} Marshall A. W., Olkin I. (1997). A new method for adding a parameter to a family of distributions with application to the exponential and Weibull families. {\sl Biometrika} 84: 641--652.

\bibitem[Mudholkar and Hutson(2000)]{Mud00} Mudholkar, G. S. and Hutson, A. D. (2000). The epsilon-skew-normal distribution for analyzing near-normal data. {\sl Journal of Statistical Planning and Inference} 83: 291--309.

\bibitem[Steutel et al.(1979)]{Steu79} Steutel, F. W., Kent, J. T., Bondesson, L. and Barndorff-Nielsen, Ole (1979). Infinite Divisibility in Theory and Practice. {\sl Scandinavian Journal of Statistics} 6: 57--64.

\bibitem[Steutel and Van Harn(2004)]{Steu04} Steutel, F. W. and Van Harn, K. (2004). {\sl Infinite Divisibility of Probability Distributions on the Real Line}. Marcel Dekker Inc., New York.

\bibitem[Wang et al.(2004)]{Wang04} Wang, J., Boyer, J. and Genton M. C. (2004). A skew symmetric representation of multivariate distributions. {\sl Statistica Sinica} 14: 1259--1270.

\end{thebibliography}
\end{document}